\pdfoutput=1
\documentclass{article}

\usepackage{amsfonts, amsmath, amssymb, amsthm, bm, booktabs, cancel, color, diagbox,
            float, graphicx, latexsym, makecell, pdfsync, setspace, url, yfonts}
            \usepackage{latexsym}
            \usepackage{authblk}

\restylefloat{table}

\theoremstyle{plain} \numberwithin{equation}{section}
\newtheorem{theorem}{Theorem}[section]

\newtheorem{corollary}[theorem]{Corollary}
\theoremstyle{dgdef}

\newcommand{\pg}{\mathrm{PG}}
\newcommand{\pgl}{\mathrm{PGL}}
\newcommand{\pgo}{\mathrm{PGO}}
\newcommand{\pGl}{\mathrm{P\Gamma L}} 
\newcommand{\nsq}{\cancel{\square}}

\newcommand{\rmQ}{\mathrm{Q}}

\newcommand{\bbF}{\mathbb{F}}

\newcommand{\bbS}{\mathbb{S}}

\newcommand{\cB}{\mathcal{B}}
\newcommand{\cC}{\mathcal{C}}

\newcommand{\cE}{\mathcal{E}}
\newcommand{\cF}{\mathcal{F}}

\newcommand{\cI}{\mathcal{I}}
\newcommand{\cK}{\mathcal{K}}

\newcommand{\cO}{\mathcal{O}}

\newcommand{\cS}{\mathcal{S}}

\newcommand{\cU}{\mathcal{U}}
\newcommand{\cV}{\mathcal{V}}
\newcommand{\cW}{\mathcal{W}}

\newcommand{\Aut}{\mathrm{Aut}}

\newcommand\blfootnote[1]{%
  \begingroup
  \renewcommand\thefootnote{}\footnote{#1}%
  \addtocounter{footnote}{-1}%
  \endgroup
}

\begin{document}
  \bibliographystyle{plain-annote}
  \title{Classification of 8-dimensional rank two commutative semifields
          \blfootnote{
            The authors acknowledge funding from the research project ``Finite Geometry with Applications
            in Algebra and Combinatorics'', funded by the Dipartimento di Tecnica e Gestione dei Sistemi Industriali
            of the Universit\'{a} di Padova.
          }
  }
  \author{Michel Lavrauw}
  \author{Morgan Rodgers}
  \affil{Universit\'a di Padova\\
    Dipartimento di tecnica e
    gestione dei sistemi industriali\\
    michel.lavrauw@unipd.it\\
    morgan.joaquin@gmail.com }
  \renewcommand\Authands{ and }

\maketitle

\begin{abstract}
  We classify the rank two commutative semifields which are 8-dimensional over their center $\bbF_{q}$.
  This is done using computational methods utilizing the connection to linear sets in $\pg(2,q^{4})$. We then
  apply our methods to complete the classification of rank two commutative semifields which are 10-dimensional
  over $\bbF_{3}$.  The implications of these results are detailed for other geometric structures such as
  semifield flocks, ovoids of parabolic quadrics, and eggs.
\end{abstract}

\section{Introduction and motivation}\label{sec:intro}

A \textit{semifield} is a possibly non-associative algebra with a one and without zero divisors.
Finite semifields are well studied objects in combinatorics and finite geometry and have
many connections to other interesting geometric structures.
They play a central role in the theory of projective planes (\cite{HuPi1973}),
generalised quadrangles (\cite{PaTh1984}), and polar spaces (\cite{Veldkamp1959}),
and have applications to perfect nonlinear functions and cryptography (\cite{BlNy2015}),
and maximum rank distance codes (\cite{Delsarte1978}).
We refer to the chapter~\cite{LaPo2011} and the references contained therein for background, definitions
and more details about these connections.

Of particular interest are commutative semifields which are of rank two over their middle nucleus,
so-called \textit{Rank Two Commutative Semifields} (R2CS)
(see~\cite{CoGa1982},~\cite{BaLa2002},~\cite[Section 5]{LaPo2011}).
The property of being commutative implies that these semifields have applications to perfect nonlinear functions
(see e.g.~\cite{BuHe2010} for a survey on planar functions and commutative semifields and for further references).
Moreover, R2CS are equivalent to semifield flocks of a quadratic cone in a 3-dimensional projective space.
We refer to the introduction of~\cite{LaPe2003} for an excellent historical overview of
the theory of flocks in finite geometry.
Consequently, R2CS are also equivalent to translation ovoids of $Q(4,q)$,
the parabolic quadric in 4-dimensional projective space.
We refer
to~\cite{FiTh1979},~\cite{Lunardon1997},~\cite{Lavrauw2005a},~\cite{Lavrauw2001},~\cite{BaBr2004},~\cite{Lavrauw2005b}
for further details on these connections.
Another, rather remarkable, application of R2CS concerns the theory of eggs and translation generalized quadrangles,
see~\cite[Section 8.7]{PaTh1984}.
As of now the only known examples of eggs in $\pg(4n-1,q)$ are either ``elementary'',
i.e.\ obtained from an oval or an ovoid by applying the technique of field reduction (\cite{LaVa2015}),
or they are obtained from a R2CS (up to dualising) (see e.g.~\cite[Chapter 3]{Lavrauw2001},~\cite{Lavrauw2005a}).

In this paper we present a computational classification of 8-dimensional
rank two commutative finite semifields (that is, 8-dimensional over their centre).
This classification relies on the bounds obtained in~\cite{BaBlLa2003} and~\cite{Lavrauw2006} on
the size of the centre in function of the dimension. Previous classification results
have been obtained for 2-dimensional semifields (\cite{Dickson1906}),
for 3-dimensional semifields (\cite{Menichetti1977} and~\cite{BaBlLa2003}),
for 4-dimensional rank two semifields (\cite{CaPoTr2006})
and for 6-dimensional rank two semifields with an extra assumption on
the size of one of the other nuclei (\cite{MaPoTr2011}).
Computational classification results have been obtained in~\cite{RuCoRa2009}, and~\cite{RuCoRa2012}.
For an overview and further classification results in the theory of finite semifields
we refer to~\cite[Section 1]{Lavrauw2013} and~\cite[Section 6]{LaPo2011}.

We begin in Section~\ref{sec:prelim} by establishing some basic terminology and
giving details on the known examples; we then explain the geometric
model we use to search for new examples of rank two commutative semifields.
In Section~\ref{sec:sublines} we determine which fields $\bbF_{q}$ satisfy
a necessary condition to be the centre of an 8-dimensional R2CS, and in
Sections~\ref{sec:planes} and~\ref{sec:higherrank} we give the results
of our exhaustive search for new examples for the field orders which meet this
necessary condition. Finally in Section~\ref{sec:implications} we give the corresponding
existence results for semifield flocks in $\pg(3,q^{4})$, translation ovoids of $Q(4,q^{4})$,
and eggs in $\pg(15,q)$.

\section{Preliminaries}\label{sec:prelim}
We use the notation and terminology from~\cite{LaPo2011}.
Given a finite semifield $\bbS$ with multiplication $(x,y)\mapsto x\circ y$,
it is natural to consider the following substructures.
The \textit{left nucleus} ${\mathbb{N}}(\bbS)$ is the set of elements $x\in \bbS$ such that
for all $y,z \in \bbS$: $x\circ(y\circ z)=(x\circ y)\circ z$.
Analogously, one defines the middle and right nucleus.
The intersection of these three subsets of $\bbS$ is called the nucleus of $\bbS$
and the intersection of the nucleus of $\bbS$ with the commutative centre of $\bbS$
is called the \textit{centre} of $\bbS$. When we mention the \textit{dimension} of a semifield,
we are referring to the dimension over its centre.
Restricting the addition and multiplication to any of these substructures one obtains a field.
The \textit{rank} of $\bbS$ is the dimension of $\bbS$ as a vector space over its middle nucleus.
Hence a \textit{rank two commutative semifield} (R2CS) is a semifield with commutative multiplication
and which is a two-dimensional vector space over its middle nucleus. A semifield $\bbS$ is commutative if
and only if $\bbS^d=\bbS$, and a semifield is called \textit{symplectic} if and only if $[\bbS^t]=[\bbS]$.

Semifields are studied up to the isotopism and their Knuth orbit.
Two semifields $\bbS_1$ and $\bbS_2$ are \textit{isotopic} if there exist non-singular linear maps
$F$, $G$, $H$ from $\bbS_1$ to $\bbS_2$ such that
$x^F \circ_2 y^G={(x\circ_1 y)}^H$ for all $x,y \in \bbS_1$.
The isotopism class of $\bbS$ is denoted by $[\bbS]$. The \textit{Knuth orbit} of a semifield $\bbS$ is
a set of at most six isotopism classes
${\mathcal{K}}({\bbS})=\{[\bbS], [\bbS^t], [\bbS^d],[\bbS^{td}],[\bbS^{dt}],[\bbS^{tdt}]\}$,
where the operations $t$ and $d$ denote the
\textit{transpose} and \textit{dual} operations obtained from the action of the transpositions in
the symmetric group $S_3$ on the indices of the cubical array of structure constants of the semifield.

To fix notation when considering a R2CS $\bbS$, we will denote the centre by $\bbF_q$,
the finite field with $q$ elements, and we will assume the left nucleus is $\bbF_{q^n}$.
This makes $\bbS$ into a $2n$-dimensional R2CS of size $q^{2n}$.


There are only three known examples of R2CS.\@(Note that by the above definition a finite field
is not a R2CS, but in some papers the finite field is also considered as an R2CS.)
We give a representation of the corresponding multiplications as binary operations
defined on $\bbF_{q^n}\times \bbF_{q^n}$ and denoted by juxtaposition $\circ$.
Also note that $n$ is necessarily at least 2 since for $n=1$ one obtains a 2-dimensional semifield
which, by Dickson~\cite{Dickson1906}, is a field.

The first example goes back to a construction by Dickson in~\cite{Dickson1906}
and exists for each odd prime power $q$ and $n\geq 2$:
\begin{eqnarray}\label{eqn:Dickson}
(x,y)\circ(u,v)=(xv+yu,yv+mx^\sigma u^\sigma),
\end{eqnarray}
where $\sigma\in \Aut(\bbF_{q^n}/\bbF_q)$ and $m$ is a non-square in $\bbF_{q^n}$.

The second family of R2CS was constructed by Cohen and Ganley~\cite{CoGa1982} and exists for $q=3$ and $n\geq 2$:
\begin{eqnarray}\label{eqn:CohenGanley}
(x,y)\circ(u,v)=(xv+yu+x^3u^3,yv+\eta x^9 u^9+\eta^{-1}xu),
\end{eqnarray}
where $\eta$ is a non-square in $\bbF_{3^n}$.
The third family is the example found by Penttila and Williams~\cite{PeWi2000} for $q=3$ and $n=5$
and has multiplication:
\begin{eqnarray}\label{eqn:PenttilaWilliams}
(x,y)\circ(u,v)=(xv+yu+x^{27}u^{27},yv+x^9 u^9).
\end{eqnarray}

Cohen and Ganley~\cite{CoGa1982} showed that R2CS in even characteristic don't exist
(again note that with our definition the finite field is not an R2CS) and that any R2CS
in odd characteristic arises from what we will
refer to as a \textit{Cohen--Ganley pair} of functions $(f,g)$: a pair of $\bbF_q$-linear functions
satisfying the property that $g^2(t)-4tf(t)$ is a non-square for all nonzero $t\in \bbF_{q^n}$, $q$ odd.
Each Cohen--Ganley pair of functions $(f,g)$ gives rise to a semifield $\bbS(f,g)$ with multiplication
\begin{eqnarray}
(x,y)\circ (u,v)=(xv+yu+f(xu),yv+g(xu)).
\end{eqnarray}
The condition that $g^2(t)-4tf(t)$ is a non-square for all nonzero $t\in \bbF_{q^n}$
is equivalent to the existence of an $\bbF_q$-linear set $\cW$ of rank $n$ whose points
have coordinates $(t,f(t),g(t))$, $t\in \bbF_{q^n}^*$, contained in the set of internal points
of the conic with equation $X_{2}^2-4X_{0}X_{1}=0$.

If $\cW$ is contained in a line then the R2CS is of Dickson type.
So we are interested in examples where $\cW$ contains an $\bbF_{q}$-subplane of $\pg(2,q^{n})$.
Using a computer search,
we complete the classification of $8$-dimensional R2CSs. This is equivalent to classifying the
semifield flocks in $\pg(3,q^{4})$ having kernel containing $\bbF_{q}$,
the translation ovoids of $Q(4,q^{4})$ with kernel containing $\bbF_{q}$, and
good eggs in $\pg(15,q)$.
We also classify the 10-dimensional R2CS with centre $\bbF_{3}$, the semifield flocks in $\pg(3,3^{5})$
with kernel $\bbF_{3}$, the translation ovoids in $Q(4,3^{5})$ with kernel $\bbF_{3}$, and the
good eggs of $\pg(19,3)$. These applications are detailed in Section~\ref{sec:implications}.

Our work relies on bounds given on the size of the centre, as a function of the dimension,
that were first given in~\cite{BaBlLa2003} and later improved in~\cite{Lavrauw2006}
by showing that in order for an $\bbF_{q}$-subplane contained in $\cI(\cC)$ to exist,
there must be an $\bbF_{q}$-subline contained in an external line of $\cC$
and made up entirely of points of $\cI(\cC)$.
\begin{theorem}[\cite{Lavrauw2006}]\label{thm:boundq}
  There are no $\bbF_{q}$-sublines contained in $\ell \cap \cI(\cC)$, where $\cC$ is a conic
  in $\pg(2,q^{n})$ and $\ell$ is an external line to $\cC$, for
  \[
    q \geq 4n^{2}-8n+2,
  \]
  and for
  \[
      q>2n^{2}-(4-2\sqrt{3})n+(3-2\sqrt{3})
  \]
  when $q$ is prime.
\end{theorem}
\begin{corollary}[\cite{Lavrauw2006}]
  No $\bbF_{q}$-subplane contained in $\cI(\cC)$ exists, where $\cC$ is a conic
  in $\pg(2,q^{n})$, for
  \[
    q \geq 4n^{2}-8n+2,
  \]
  and for
  \[
      q>2n^{2}-(4-2\sqrt{3})n+(3-2\sqrt{3})
  \]
  when $q$ is prime.
\end{corollary}

Let $q$ be odd, and consider the conic $\cC$ in $\pg(2,q^{n})$ defined by the quadratic form
$\rmQ: X_{0}X_{1} - X_{2}^{2}$.  Notice that the point $(0,0,1)$ lies on the tangent $[1,0,0]$,
so this point is external. Since $\rmQ(0,0,1) = -1$, we have that the internal points $\cI(\cC)$
are those for which $-\rmQ(\bm{v}) \in \nsq$.
The stabilizer $G = \pgo(3,q^{n})$ of $\cC$ in $\pgl(3,q^{n})$
has order $q^{n}(q^{2n}-1)$, and contains all matrices of the form
\[
\begin{bmatrix} a^{2} & b^{2} & ab\\ c^{2} & d^{2} & cd\\ 2ac & abd & ad+bc \end{bmatrix}
\]
where $ad-bc \neq 0$ (vector multiplication is from the left).

We have the following, due to Payne~\cite{Payne}:
\begin{theorem}\leavevmode
  \begin{enumerate}
  \item $G$ is sharply triply transitive on the points of $\cC$;
  \item $G$ is transitive on $\cI(\cC)$;
  \item $G$ is transitive on $\cE(\cC)$;
  \item $G$ is sharply triply transitive on the tangent lines to $\cC$;
  \item $G$ is transitive on the external lines to $\cC$;
  \item $G$ is transitive on the secant lines to $\cC$;
  \item $G$ is transitive on the point-line pairs $(\bm{p}, \ell)$, where $\bm{p}$ is an
    external point on the exterior (resp., secant) line $\ell$.
    The subgroup of $G$ fixing such a pair has order $4$.
  \item $G$ is transitive on the point-line pairs $(\bm{p}, \ell)$, where $\bm{p}$ is an
    internal point on the exterior (resp., secant) line $\ell$.
    The subgroup of $G$ fixing such a pair has order $4$.
  \end{enumerate}
\end{theorem}

\section{Existence of sublines contained in $\cI(\cC)$}\label{sec:sublines}
Our first goal is to determine precisely the values of $q$
for which there exists an $\bbF_{q}$-subline
contained in an external line to a conic $\cC$ in $\pg(2,q^{n})$
consisting entirely of internal points of $\cC$.
To accomplish this we choose $\eta$ so that $-\eta \in \nsq$ and
$-\eta-1 \in \square$; then $\bm{x} = (1, \eta, 0) \in \cI(\cC)$
and $\ell_{e} = \langle (1, \eta, 0),\ (0, -2 \eta, 1) \rangle$
is an external line.
The stabilizer in $G$ of $\bm{x}$ has order $2(q^{n}+1)$,
and contains the following (normalized) matrices:
\[
G_{\bm{x}} =
  \left\{    
    \begin{bmatrix}
      1 & 0 &     0 \\
      0 & 1 &     0 \\
      0 & 0 & \pm 1
    \end{bmatrix}
    \right\} \cup 
    \left\{ 
    \begin{bmatrix}
      a^{2} & 1 & a \\
      \frac{1}{\eta^{2}} & a^{2} & -\frac{a}{\eta} \\
      \pm \frac{2a}{\eta} & \mp 2a & \pm \left( \frac{1}{\eta} - a^{2} \right)
    \end{bmatrix} : a \in \mathbb{E}
  \right\}    
\]

Now considering the external line $\ell_{e} = \langle (1, \eta, 0),\ (0, -2 \eta, 1) \rangle$
on $\bm{x}$, the subgroup of $G$ stabilizing both $\bm{x}$ and $\ell_{e}$ has order 4 and is given by
\[
G_{\bm{x}, \ell_{e}} =
  \left\langle    
    \begin{bmatrix}
      0                  & 1 & 0 \\
      \frac{1}{\eta^{2}} & 0 & 0 \\
      0                  & 0 & -\frac{1}{\eta}
    \end{bmatrix},\
    \begin{bmatrix}
      1                  &  1 & 1 \\
      \frac{1}{\eta^{2}} &  1 & -\frac{1}{\eta} \\
      \frac{2}{\eta}     & -2 & \left( \frac{1}{\eta} - 1 \right)
    \end{bmatrix}
  \right\rangle.    
\]
Note that the second generator given for this group fixes $\ell_{e}$ pointwise.

Since $G$ acts transitively on pairs $(\bm{p}, \ell)$, where $\bm{p}$ is an internal point on an external line $\ell$,
it is sufficient to look for sublines of $\ell_{e}$ containing $\bm{x}$ and contained in $\cI(\cC)$.
A subline is determined by three collinear points; so a subline of $\ell_{e}$ contained in $\cI(\cC)$ is determined by
$\bm{x}$, $\bm{y}$, and $\bm{x}+ \mu \bm{y}$ for some $\bm{y} \in (\ell_{e} \cap \cI(\cC)) \setminus \{\bm{x}\}$
and $\mu \in \bbF_{q^{n}}^{*}$ satisfying $-\rmQ(\bm{x} + \lambda \mu \bm{y}) \in \nsq$ for all $\lambda \in \bbF_{q}$.
The subline determined by these three points is given by
$\{ \bm{y} \} \cup \{ x + \lambda \mu \bm{y} \ : \ \lambda \in \bbF_{q}\}$.

Using the basis $\{ \bm{v}_{1} = (1, \eta, 0),\ \bm{v}_{2}=(0, -2 \eta, 1) \}$,
we can associate $\ell_{e}$ with $\pg(1,q^{n})$ having the induced quadratic form
$\rmQ_{\ell_{e}}(x_{1}\bm{v}_{1}+x_{2}\bm{v}_{2}) = \eta x_{1}^{2} - 2\eta x_{1}x_{2} - x_{2}^{2}$
(this form is anisotropic, and is just used to separate the points of $\ell_{e}$ into $\cI(\cC)$ and $\cE(\cC)$).

We want to take a point in $(\ell_{e} \cap \cI(\cC))\setminus \{\bm{x} \}$.
Since $- \rmQ(\bm{v}_{2}) = 1 \in \square$, $\bm{v_{2}} \not\in \cI(\cC)$;
therefore we will define $\bm{y}_{b} = \bm{v}_{1} + b \bm{v}_{2}$ with $b \neq 0$.
Now we will have $\bm{y}_{b} \in \cI(\cC)$ as long as $b$ satisfies $b^{2}+ 2\eta b - \eta \in \nsq$.
Let $\cB = \{ s : s \in \mathbb{F}_{q^{n}} \mid s^{2}+ 2\eta s - \eta \in \nsq \}$, then we have that
\[
(\ell_{e} \cap \cI(\cC)) = \{\bm{x} \} \cup \{ \bm{y}_{b} : b \in \cB \}.
\]

Now for our choices of $\mu$, instead of letting $\mu$ range over all possible values in $\bbF_{q^{n}}^{*}$,
it is sufficient to consider a set $\cS$ of representatives of $\bbF_{q^{n}}^{*}/\bbF_{q}^{*}$.
Notice that
\[
  \bm{x} + \lambda \mu \bm{y}_{b} = \bm{v}_{1} + \lambda \mu (\bm{v}_{1} + b \bm{v}_{2})
    = (1+\lambda \mu)\bm{v}_{1} + \lambda \mu b \bm{v}_{2};
\]
normalizing this vector to $\bm{v}_{1} + \frac{\lambda \mu}{1+\lambda \mu}b \bm{v}_{2}$,
we see that it is contained in $\cI(\cC)$ if and only if $\frac{\lambda \mu}{1+\lambda \mu}b \in \cB$.
This is equivalent to having
\[
(2b-1)\mu^{2}\eta\lambda^{2}+2(b-1)\mu\eta\lambda +b^{2}-\eta \in \nsq
\]
for all $\lambda \in \bbF_{q}$.

To find sublines efficiently, we compute the set $\cB$ and, for each value of $\mu \in \cS$,
the sequence $[\frac{\lambda \mu}{1+\lambda \mu} : \lambda \in \bbF_{q}^{*}]$. Then for each pair
$(b,\mu) \in \cB \times \cS$,
we check whether $\frac{\lambda \mu}{1+\lambda \mu}b \in \cB$ for all $\lambda \in \bbF_{q}$.
In this way, we obtain the number of $\bbF_{q}$-sublines of $\ell_{e}$ containing $\bf{x}$ and completely contained
in $\cI(\cC)$ for $\cC$ a conic in $\pg(2,q^n)$ for $n=3$ and $n=4$.
We also obtain some partial results for $n=5$, however it was impossible for us to complete our computations
for $q \in \{37, 41, 43, 49\}$ in this case. Our results for $n=3$ agree with those found in~\cite{BlThVM98}.
Notice that from the bounds given by Theorem~\ref{thm:boundq},
when $n=3$ we only need to consider $q < 14$; for $n=4$ we only need to consider $q < 30$;
and when $n=5$ we only need to consider $q < 47$, along with $q=49$. In the table, a $0$ indicates
that no subline was found, while a dash indicates that the existence is ruled out by the theoretical bound.\\
\begin{table}[ht]
  \begin{tabular}{rrrr}
        & \multicolumn{3}{c}{Number of sublines on $\ell_{e}$}\\
    $q$ & $n=3$ & $n=4$ & $n=5$\\
    \toprule
    3   & 12 &  120 &   1200\\
    5   & 12 &  600 &  15072\\
    7   & 24 &  912 &  52080\\
    9   &  0 & 1040 &  91880\\
    11  &  0 &  744 & 115572\\
    13  &  0 &  504 & 102340\\
    17  &  - &   72 & $\geq 1$\\
    19  &  - &   80 & $\geq 1$\\
    23  &  - &    0 & $\geq 1$\\
    25  &  - &    0 & $\geq 1$\\
    27  &  - &    0 & $\geq 1$\\
    29  &  - &    0 & $\geq 1$\\
    31  &  - &    - & $\geq 1$\\
  \end{tabular}
\end{table}

Our computational results show the following.
\begin{theorem}\label{thm:sublines}
  If there exists an $\bbF_{q}$-subline in $\pg(2,q^{4})$ contained in $\ell \cap \cI(\cC)$ for some
  conic $\cC$, where $\ell$ is an external line to $\cC$, then $q \leq 19$.
\end{theorem}

\begin{table}[ht]
  \begin{tabular}{llllllll}
    \multicolumn{8}{c}{Sublines for $3^4$ ($(b, \mu)$ pairs) containing $\bm{x} = (1,\alpha,0)$}\\
    \multicolumn{8}{c}{Minimal polynomial of $\alpha$: $x^{4}+2x^{3}+2$}\\
    $b$ & $\mu$ & & $b$ & $\mu$ & & $b$ & $\mu$\\
    \toprule
    $\alpha^{2}$:   & $\alpha^{54}$, $\alpha^{56}$ & &
    $\alpha^{5}$:   & $\alpha^{27}$, $\alpha^{48}$ & &
    $\alpha^{6}$:   & $\alpha$\\
    $\alpha^{11}$:  & $\alpha^{2}$ & &
    $\alpha^{13}$:  & $\alpha^{49}$ & &
    $\alpha^{14}$:  & $\alpha^{27}$\\
    $\alpha^{15}$:  & $\alpha^{30}$ & &
    $\alpha^{16}$:  & $\alpha^{10}$ & &
    $\alpha^{17}$:  & $\alpha^{27}$\\
    $\alpha^{18}$:  & $\alpha^{27}$, $\alpha^{72}$ & &
    $\alpha^{20}$:  & $\alpha$, $\alpha^{37}$, $\alpha^{48}$ & &
    $\alpha^{21}$:  & $\alpha^{5}$, $\alpha^{58}$\\
    $\alpha^{23}$:  & $\alpha^{5}$, $\alpha^{49}$ & &
    $\alpha^{28}$:  & $\alpha^{12}$ & &
    $\alpha^{30}$:  & $\alpha^{23}$, $\alpha^{54}$\\
    $\alpha^{32}$:  & $\alpha^{10}$ & &
    $\alpha^{38}$:  & $\alpha^{56}$ & &
    $\alpha^{43}$:  & $\alpha^{30}$\\
    $\alpha^{44}$:  & $\alpha^{38}$ & &
    $\alpha^{45}$:  & $\alpha^{15}$, $\alpha^{66}$ & &
    $\alpha^{46}$:  & $\alpha^{13}$\\
    $\alpha^{47}$:  & $\alpha^{22}$ & &
    $\alpha^{49}$:  & $\alpha^{30}$ & &
    $\alpha^{50}$:  & $\alpha^{5}$\\
    $\alpha^{51}$:  & $\alpha^{23}$, $\alpha^{49}$ & &
    $\alpha^{54}$:  & $\alpha^{27}$, $\alpha^{75}$ & &
    $\alpha^{55}$:  & $\alpha^{66}$\\
    $\alpha^{57}$:  & $\alpha^{10}$ & &
    $\alpha^{58}$:  & $\alpha^{30}$ & &
    $\alpha^{60}$:  & $\alpha^{54}$\\
    $\alpha^{67}$:  & $\alpha^{56}$ & &
    $\alpha^{68}$:  & $\alpha^{2}$, $\alpha^{27}$, $\alpha^{54}$ & &
    $\alpha^{70}$:  & $\alpha^{5}$, $\alpha^{10}$, $\alpha^{73}$\\
    $\alpha^{72}$:  & $\alpha^{27}$ & &
    $\alpha^{73}$:  & $\alpha^{39}$ & &
    $\alpha^{75}$:  & $\alpha^{75}$\\
    $\alpha^{76}$:  & $\alpha^{27}$\\
  \end{tabular}
\end{table}


\section{Finding subplanes}\label{sec:planes}
Our next goal is to determine, given the existence of the necessary $\bbF_{q}$-subline, when there
exist $\bbF_{q}$-subplanes of $\pg(2,q^{n})$ contained in $\cI(\cC)$.  An $\bbF_{q}$-subplane is
completely determined by a quadrangle, so more generally, two $\bbF_{q}$-sublines
that are not contained in a common line of $\pg(2,q^{n})$ will determine an $\bbF_{q}$-subplane
of $\pg(2,q^{n})$.

To determine the existence of $\bbF_{q}$-subplanes contained in $\cI(\cC)$,
we first fix the point $\bm{x}$ and then find all of the $\bbF_{q}$-sublines containing $\bm{x}$
which are completely contained in $\cI(\cC)$ (those spanning an external line to $\cC$ as well as
those spanning a secant line).
Then, for each pair of $\bbF_{q}$-sublines through $\bm{x}$ (not contained in a common line of $\pg(2,q^{n})$),
we test whether the $\bbF_{q}$-subplane they determine is a subset of $\cI(\cC)$.

In the previous section we give details on finding the $\bbF_{q}$-sublines of the external line
$\ell_{e}$ on $\bm{x}$ which are contained in $\cI(\cC)$.
Once we have these sublines, we compute their images under $G_{x}$
to get all of the $\bbF_{q}$-sublines on $\bm{x}$ contained in $\cI(\cC)$
generating an external line to $\cC$.
We then repeat this process beginning with
the secant line $\ell_{s} = \langle (1,0,0), (0,1,0) \rangle$ on $\bm{x}$.
Since all $\bbF_{q}$-sublines are assumed to contain $\bm{x}$,
and an $\bbF_{q}$-subline is determined by $3$ points,
we save the sublines as an ordered pair $\{@ \bm{y}, \bm{z} @\}$
where $\{\bm{x},\bm{y},\bm{z}\}$ determines the subline.

The real computationally intensive aspect of our work concerns
the determination of whether two sublines form a compatible pair,
that is, if the two sublines determine a rank $3$ $\bbF_{q}$-linear set
which is contained in $\cI(\cC)$.  For two $\bbF_{q}$ sublines $\ell_{1}$ and $\ell_{2}$
generated by $\{\bm{x}, \bm{y}_{1}, \bm{z}_{1}\}$ and $\{\bm{x}, \bm{y}_{2}, \bm{z}_{2}\}$,
respectively, we first compute values $\mu_{1}$ and $\mu_{2}$ so that the $\bbF_{q}$-subplane
spanned by these two lines is given by
$\langle \bm{x}, \mu_{1}\bm{y}_{1}, \mu_{2}\bm{y}_{2} \rangle_{q}$.
Then we test that $\lambda \mu_{1}\bm{y}_{1}+\bm{y}_{2} \in \cI(\cC)$ for all $\lambda \in \bbF_{q}^{*}$,
and that $\bm{x} + \lambda_{1} \mu_{1}\bm{y}_{1}+\lambda_{2} \mu_{2} \bm{y}_{2} \in \cI(\cC)$ for all
$\lambda_{1}$, $\lambda_{2} \in \bbF_{q}^{*}$. If these conditions are satisfied,
then $\ell_{1}$ and $\ell_{2}$ generate an $\bbF_{q}$-subplane contained in $\cI(\cC)$.

Our computational work proves the following.
\begin{theorem}
  No $\bbF_{q}$-subplane contained in $\cI(\cC)$ exists,
  where $\cC$ is a conic in $\pg(2,q^4)$, unless $q=3$.
\end{theorem}
With $n=4$ and $q = 3$ we find 13 $\bbF_{3}$-subplanes (up to conjugacy in $\pGl(3,3^{4})$)
contained in $\cI(\cC)$, 10 of which can be embedded in the linear set associated with the Cohen--Ganley
semifield.

\section{Linear sets of higher rank}\label{sec:higherrank}
To put together rank $4$ $\bbF_{q}$-linear sets, we first need to find \textit{all} the rank $3$ linear sets
(not just the subplanes). It is fairly easy to find the examples that are contained in a line.
Each rank $3$ linear set is saved as an ordered pair of $\bbF_{q}$-linear lines containing $\bm{x}$.
Then, for each $\bbF_{q}$-subline contained in either $\ell_{e}$ or $\ell_{s}$, we compile the set
$\Pi_{\ell}$ of rank $3$ $\bbF_{q}$-linear sets whose first generating line is $\ell$.
We form a graph $\Gamma_{\ell}$ on $\Pi_{\ell}$, where two planes $\pi_{1}$, $\pi_{2} \in \Pi_{\ell}$
are adjacent if their second generating lines together generate a rank $3$ $\bbF_{q}$-linear set
contained in $\cI(\cC)$. Then a clique of size $q(q+1)$ in $\Gamma_{\ell}$ corresponds to a rank $4$
$\bbF_{q}$-linear set contained in $\cI(\cC)$.

Running this algorithm using the rank $3$ linear sets found in $\pg(2,3^{4})$, we find
$174$ rank $4$ $\bbF_{q}$-linear sets contained in $\cI(\cC)$ that contain an $\bbF_{q}$-subplane.
They are all equivalent up to isomorphism, corresponding to a semifield of Cohen--Ganley type.
We are also able to run this algorithm in $\pg(2,3^5)$, using an increased clique size to look for
rank $5$ $\bbF_{q}$-linear sets; here all examples found correspond to a semifield of Cohen--Ganley type
or else to the example due to Penttila and Williams.
\begin{theorem}
  An 8-dimensional R2CS is either a Dickson semifield, or of Cohen--Ganley type (with centre $\bbF_{3}$).
\end{theorem}
\begin{theorem}
  A 10-dimensional R2CS with centre $\bbF_{3}$ is either a Dickson semifield, of Cohen--Ganley type,
  or Penttila--Williams.
\end{theorem}

\section{Implications of our results}\label{sec:implications}
There are many connections between R2CS and various geometric objects.  Here we give details on some of these
connections, and state the implications of our classification of R2CS in these other settings.
\subsection{Semifield flocks}
A flock of a quadratic cone $\cK$ of $\pg(3,q^{n})$ having vertex $v$ is a partition of $\cK \setminus \{v\}$
into $q^{n}$ conics.  We let $v = (0,0,0,1)$ and let the conic $\cC$ in the plane $\pi = [0,0,0,1]$
be the base of the cone.  Then the planes of the flock can be written as
\[
  \{ \pi_{t} : tX_{0} +f(t)X_{1} + g(t)X_{2}+X_{3}=0 \mid t \in \bbF_{q^{n}} \}
\]
for some $f, g: \bbF_{q^{n}} \to \bbF_{q^{n}}$; we denote such a flock by $\cF(f,g)$.

A flock corresponds to a set
\[
  \cW = \{ (t, f(t), g(t))  \mid t \in \bbF_{q^{n}} \}
\]
 of interior points of a conic $\cC^{\prime}$ in $\pg(2,q^{n})$ (see~\cite{Thas1987}).
 If $f$ and $g$ are linear over a subfield
 of $\bbF_{q^{n}}$ (i.e.\ if $(f,g)$ is a Cohen--Ganley pair) then we say $\cF(f,g)$ is a \textit{semifield} flock.
 The maximal subfield of $\bbF_{q^{n}}$ for which $f$ and $g$ are linear is called the \textit{kernel}
 of the semifield flock.  Notice that if $\cF(f,g)$ is a semifield flock of $\pg(3,q^{n})$ with kernel $\bbF_{q}$
 then $\cW$ is a rank $n$ $\bbF_{q}$-linear set contained in the set of interior points of a conic in $\pg(2,q^{n})$,
 so such a semifield flock is equivalent to a $2n$-dimensional R2CS with centre $\bbF_{q}$.
 Furthermore if $\cW$ is contained in a line then $\cF$ is of Kantor--Knuth type~\cite{Thas1987};
 this corresponds to a R2CS of Dickson type.

 \begin{corollary}
   A semifield flock of $\pg(3,q^{4})$ with kernel $\bbF_{q}$ is of Kantor--Knuth type
   or of Cohen--Ganley type (with kernel $\bbF_{3}$).
 \end{corollary}
 \begin{corollary}
   A semifield flock of $\pg(3,3^{5})$ with kernel $\bbF_{3}$ is of Kantor--Knuth type, of Cohen--Ganley type,
   or of Penttila--Williams type.
 \end{corollary}

\subsection{Ovoids of the parabolic quadric in $4$-dimensional projective space}
The parabolic quadric $Q(4,s)$ is the incidence structure of points and lines of a nondegenerate quadric
in $\pg(4,s)$. The quadric $Q(4,s)$ is also
an example of a generalized quadrangle, and is known to be isomorphic to
the example $T_{2}(\cC)$ constructed from a conic $\cC$ in $\pg(2,s)$, see~\cite{PaTh1984}.

A set of $s^{2}+1$ points $\cO$ of $Q(4,s)$ is called an \textit{ovoid} if no
two points of $\cO$ are collinear in $Q(4,s)$.  An ovoid $\cO$ in $Q(4,s)$ is a \textit{translation ovoid}
if there is a point $\bm{p} \in \cO$ and a group $G$ of collineations of $Q(4,s)$ stabilizing $\cO$,
fixing $\bm{p}$, and acting regularly on the points of $\cO\setminus \{\bm{p}\}$.  This group $G$ is necessarily
elementary abelian, and hence is a vector space over some subfield $\bbF_{q}$ of $\bbF_{s}$; the largest such
subfield is called the \textit{kernel} of the translation ovoid. Put $n = [\bbF_{s} : \bbF_{q}]$, so $s = q^{n}$.

By~\cite[Section 3.2]{Lavrauw2006}, the classification result from Theorem~\ref{thm:sublines}
has the following applications to ovoids of $Q(4,q^{4})$.  Given an ovoid $\cO$ of $Q(4,q^{n})$,
for each point $\bm{p} \in \cO$,
fix some conic $\cC_{\bm{p}}$ contained in the cone $\bm{p}^{\perp} \cap Q(4,q^{n})$;
we will denote the plane containing $\cC_{\bm{p}}$ by $\pi_{\bm{p}}$.
Then we can consider $Q(4,q^{n}) \simeq T_{2}(\cC_{\bm{p}})$.
In this model, $\bm{p}$ corresponds to the point $(\infty)$, and the points of $\cO \setminus \{\bm{p}\}$ correspond
to a set $\cV_{\bm{p}}$ of $q^{2n}$ affine points. Each two points of $\cV_{\bm{p}}$ span a line
intersecting the plane $\pi_{\bm{p}}$ in a point not on $\cC_{\bm{p}}$. Define
\[
\cU_{\bm{p}} = \{ \langle x,y \rangle \cap \pi_{\bm{p}} \mid x,y \in \cV_{\bm{p}}\}.
\]
If the set $\cU_{\bm{p}}$ contains a dual $\bbF_{q}$-subline on an internal point with respect to $\cC_{\bm{p}}$,
then dualising over $\bbF_{q}$, we have an $\bbF_{q}$-subline spanning an external line with respect to $\cC_{\bm{p}}$.
This gives us the following.
\begin{corollary}
  If $\cO$ is an ovoid in $Q(4,q^4)$, $q$ odd, and $\cU_{\bm{p}}$ contains a dual $\bbF_{q}$-subline
  on an internal point of $\cC_{\bm{p}}$, for some point $\bm{p} \in \cO$, then
  $q \leq 19$.
\end{corollary}
If $\cO$ is a translation ovoid of $Q(4,q^{n})$ with respect to the point $\bm{p}$ having kernel $\bbF_{q}$,
then the set $\cU_{\bm{p}}$ is a rank $2n$ $\bbF_{q}$-linear set, and its dual is a rank $n$ $\bbF_{q}$-linear
set contained in $\cI(\cC_{\bm{p}})$.
\begin{corollary}
  A translation ovoid in $Q(4,q^{4})$ with kernel $\bbF_{q}$ is either a Kantor ovoid,
  or a Thas--Payne ovoid (with $q=3$).
\end{corollary}
\begin{corollary}
  A translation ovoid in $Q(4,3^5)$ with kernel $\bbF_{3}$ is either a Kantor ovoid,
  a Thas--Payne ovoid, or Penttila--Williams ovoid.
\end{corollary}


\subsection{Eggs}
We define an \textit{egg} $\cE$ in $\pg(4n-1,q)$ to be a partial $(n-1)$-spread of size $q^{2n}+1$
such that every 3 elements of $\cE$ span a $(3n-1)$-space and, for every element $E\in \cE$, there exists
a $(3n-1)$-space, denoted $T_{E}$ and called the \textit{tangent space of $\cE$ at $E$}, containing
$E$ and disjoint from every other egg element.  An egg is called \textit{good at an element $E \in \cE$}
if every $(3n-1)$-space containing $E$ and at least two other elements of $\cE$ contains exactly
$q^{n}+1$ elements of $\cE$. We say that an egg $\cE$ of $\pg(4n-1,q)$ is a \textit{good egg}
if there exists an element $E \in \cE$ for which $\cE$ is good at $E$.  The standard example of
an egg in $\pg(4n-1,q)$ is obtained by applying field reduction to an ovoid of $\pg(3,q^{n})$;
an egg that can be obtained in this way is called \textit{elementary}.

It is shown in~\cite{Thas1999} (see~\cite{LaPe2001} for a shorter direct proof) that
good eggs of $\pg(4n-1,q)$, $q$ odd, are equivalent to
semifield flocks of $\pg(3,q^{n})$ with kernel containing $\mathbb{F}_{q}$. This gives us the
following result.
\begin{corollary}
  If $\cE$ is a good egg of $\pg(15,q)$ with kernel $\bbF_{q}$, $q$ odd,
  then $\cE$ is either elementary, of Kantor--Knuth type, or Cohen--Ganley.
\end{corollary}
Even if we do not assume that the egg $\cE$ has a good element, it is shown in~\cite{Lavrauw2006}
that an egg with certain properties implies the existence of an $\bbF_{q}$-subline contained
in the set of interior points of a conic $\cC$ in $\pg(2,q^{n})$ which spans an external
line with respect to $\cC$, giving the following result.
\begin{corollary}
Let $\cE$ be an egg of $\pg(15,q)$, $q$ odd. If there exists an $11$-space $\rho$
containing an elementary pseudo-oval $\cO_{q}$ contained in $\cE$ corresponding to
a conic $\cC$ of $\pg(2,q^{4})$, and there is a tangent space intersecting $\rho$
in a $7$-space $\cU$ whose associated $\bbF_{q}$-linear set in
$\langle \cC \rangle \simeq \pg(2,q^{n})$ contains a dual $\bbF_{q}$-subline
on an internal point w.r.t.\ $\cC$, then $q \leq 19$.
\end{corollary}
\begin{corollary}
  If $\cE$ is a good egg of $\pg(19,3)$ with kernel $\bbF_{3}$,
  then $\cE$ is either elementary, Kantor--Knuth,  Cohen--Ganley, or Penttila--Williams.
\end{corollary}

\bibliography{R2CSarXiv}
\end{document}